\documentclass[a4paper, reqno]{amsart}

\usepackage{amssymb}
\usepackage{amsmath}
\usepackage{amsthm}
\usepackage{mathtools}
\usepackage{complexity}
\usepackage{tikz}
\usetikzlibrary{positioning}
\usetikzlibrary{decorations.pathreplacing}
\usetikzlibrary{arrows.meta}
\usepackage{hyperref}
\usepackage{enumerate}
\usepackage{cleveref}

\DeclareMathOperator{\GL}{GL}

\theoremstyle{plain}
\newtheorem{theorem}{Theorem}[section]

\newtheorem{lemma}[theorem]{Lemma}

\newtheorem{fact}[theorem]{Fact}
\newtheorem{remark}[theorem]{Remark}

\renewcommand{\poly}{\mathrm{poly}}
\renewcommand{\C}{\mathbb{C}}

\newcommand{\BS}{\mathrm{BS}}
\newcommand{\Z}{\mathbb{Z}}
\newcommand{\N}{\mathbb{N}}

\newcommand{\shift}[2]{S_{#2}}

\newcommand{\PowWP}{\mathrm{PowerWP}}
\newcommand{\KP}{\mathrm{Knapsack}}
\newcommand{\gmat}[2]{\begin{pmatrix} #2 & #1 \\ 0 & 1 \end{pmatrix}}
\newcommand{\gsmat}[2]{\begin{psmallmatrix} #2 & #1 \\ 0 & 1 \end{psmallmatrix}}
\newcommand{\Sol}{\mathrm{sol}}
\newcommand{\ExpEq}{\mathrm{ExpEq}}

\renewcommand{\deg}{\mathrm{deg}}
\newcommand{\height}{\mathrm{height}}

\begin{document}

\title[Knapsack and power word problem in Baumslag-Solitar groups]{Knapsack and the power word problem in solvable Baumslag-Solitar groups}

\author[M.~Ganardi]{Moses Ganardi}
\email{ganardi@mpi-sws.org}
\author[M.~Lohrey]{Markus Lohrey}
\email{lohrey@eti.uni-siegen.de}
\author[G.~Zetzsche]{Georg Zetzsche}
\email{georg@mpi-sws.org}

\address[Moses Ganardi, Georg~Zetzsche]{Max Planck Institute for Software Systems, Kaiserslautern, Germany}
\address[Markus~Lohrey]{Universit{\"a}t Siegen, Germany}

\begin{abstract}
	We prove that the power word problem for certain metabelian subgroups of $\GL(2,\mathbb{C})$ (including the 
	solvable Baumslag-Solitar groups $\BS(1,q) = \langle a,t \mid t a t^{-1} = a^q \rangle$) belongs to the circuit complexity class
	$\TC^0$. In the power word problem, the input consists of group elements $g_1, \ldots, g_d$ and binary encoded integers
	$n_1, \ldots, n_d$ and it is asked whether $g_1^{n_1} \cdots g_d^{n_d} = 1$ holds.	Moreover, we prove that the knapsack problem
	for $\BS(1,q)$ is $\NP$-complete. In the knapsack problem, the input consists of group elements $g_1, \ldots, g_d,h$ 
	and it is asked whether the equation $g_1^{x_1} \cdots g_d^{x_d} = h$ has a solution in $\N^d$. For the more general case of a system
	of so-called exponent equations, where the exponent variables $x_i$ can occur multiple times, we show that solvability is undecidable for 
	$\BS(1,q)$.
	\end{abstract}



\maketitle
	
\section{Introduction}

\subsection{The power word problem} The study of multiplicative identities and equations has a long tradition in computational algebra, and has recently
been extended to the non-abelian case. Here, the multiplicative identities we have in mind have the form $g_1^{n_1} g_2^{n_2} \cdots g_d^{n_d} = 1$,
where $g_1, \ldots, g_d $ are elements of a group $G$ and $n_1, n_2, \ldots, n_d \in \N$ are non-negative integers (we may also allow negative
$n_i$, but this makes no difference, since we can replace a $g_i$ by its inverse $g_i^{-1}$). Typically, the numbers $n_i$ are given in binary
representation, whereas the representation of the group elements $g_i$ depends on the underlying group $G$. Here, we consider the case where
$G$ is a finitely generated (f.g.~for short) group, and elements of $G$ are represented by finite words over a fixed generating set $\Sigma$ (the concrete choice 
of $\Sigma$ is not relevant). In this setting, the question whether $g_1^{n_1} g_2^{n_2} \cdots g_d^{n_d} = 1$ is a true identity has been recently
introduced as the {\em power word problem} for $G$ \cite{LoWe19}. It extends the classical word problem for $G$ (does a 
given word over the group generators represent the group identity?) in the sense that the word problem trivially reduces to the power word
problem (take an identity $w^1 = 1$). Recent results on the power word problem in specific f.g.~groups are:
\begin{itemize}
\item For every f.g.~free group the power word problem belongs to deterministic logspace \cite{LoWe19}. This result has been recently generalized in \cite{StoberW22}, where
it is shown that the power word problem in a fixed graph product of groups is logspace-reducible (even $\mathsf{AC}^0$-Turing-reducible) to the word problem of the free group
of rank two and the power word problem of the base groups of the graph product.
\item For the following groups the power word problem belongs to the circuit complexity class $\TC^0$:\footnote{$\TC^0$ is a very small complexity
class within polynomial time; see Section~\ref{sec-complexity} for more details. In this paper, $\TC^0$ always
refers to the $\mathsf{DLOGTIME}$-uniform version.}  f.g.~nilpotent groups \cite{LoWe19}, iterated
wreath products of f.g.~free abelian groups and (as a consequence of the latter) free solvable groups~\cite{FGLZ20}.
\item If $G$ is a so-called uniformly efficiently non-solvable group (this is a large class of non-solvable groups that was recently introduced
in \cite{BFLW19} and that includes all finite non-solvable groups and f.g.~free non-abelian groups) then the power word problem for the wreath product $G \wr \mathbb{Z}$
is $\coNP$-hard~\cite{FGLZ20}. As a consequence, the power word problem for Thompson's group $F$ is $\coNP$-complete~\cite{FGLZ20}.
\end{itemize}
Historically, the power word problem appeared earlier in the area of computational (commutative) algebra. Ge \cite{Ge93} 
proved that one can check in polynomial time whether an identity
$\alpha_1^{n_1} \alpha_2^{n_2} \cdots \alpha_d^{n_d} = 1$, where the $n_i$ are binary encoded integers and the $\alpha_i$ are from an algebraic number field
(and suitable encoded), holds. 

In this paper we investigate the power word problem for certain 2-generated subgroups of $\GL(2,\mathbb{C})$: for a fixed complex number $\alpha \in \mathbb{C} \setminus \{0\}$
we consider the group $T(\alpha)$ generated by the two matrices
\[
   \begin{pmatrix} 1 & 1 \\ 0 & 1 \end{pmatrix}  \quad\text{and}\quad   \begin{pmatrix} \alpha & 0 \\ 0 & 1 \end{pmatrix}.
\]
In a purely group theoretic context, these groups were studied in \cite{Guyot12,Guyot18}.
Important special cases are the wreath product $\mathbb{Z} \wr \mathbb{Z}$ (for $\alpha$ transcendental) and 
the solvable Baumslag-Solitar groups $\BS(1,q) = \langle a,t \mid t a t^{-1} = a^q \rangle$ (for $\alpha = q \geq 2$ an integer).
Our first main result states that the power word problem for every group
$T(\alpha)$  belongs to $\TC^0$ (Theorem~\ref{thm-powerWP}). For the word problem of $T(\alpha)$, membership in $\TC^0$ follows from \cite{KonigL15}\footnote{For the special case $\BS(1,q)$ membership of the word problem in $\TC^0$ was shown in \cite{Weiss15b}.}
since $T(\alpha)$ is a linear solvable group. 
Theorem~\ref{thm-powerWP} is directly related to recent results on matrix powering problems \cite{AllenderBD14,GalbyOW15}. These problems
can be quite difficult to analyze. For instance, it is not known whether a certain bit of the $(1,1)$-entry of a matrix power $A^n$ 
can be computed in polynomial time, when $n$ is given in binary notation and $A$ is a $(2\times 2)$-matrix over $\mathbb{Z}$. The related problem of checking whether the $(1,1)$-entry (or any other entry) of $A^n$ is positive can be solved in polynomial time by \cite{GalbyOW15}.

\subsection{The knapsack problem} If one replaces in the power word problem the exponents $n_i$ by pairwise different variables $x_i$ and the right-hand side $1$ by
an arbitrary group element $h \in G$, one obtains a so-called knapsack equation $g_1^{x_1} g_2^{x_2} \cdots g_d^{x_d} = h$.
The question, whether such an equation has a solution in $\N^d$ is known as the {\em knapsack problem} for $G$. 
In the general context of finitely generated groups the knapsack problem has been introduced by Myasnikov, Nikolaev, and Ushakov  \cite{MyNiUs14}. 
 As for the power word problem, this problem has been studied in the commutative setting before. For the case $G = \Z$ one obtains a variant of
 the classical $\NP$-complete knapsack problem; a proof of the $\NP$-hardness of our variant of the knapsack problem for the integers
 can be found in \cite{Haa11}. For this hardness result it is important that integers are represented in binary notation. For unary encoded
 integers the complexity of the knapsack problem goes down to $\TC^0$. For the case that the $g_i$ are commuting matrices over an algebraic
 number field, the knapsack problem has been studied in \cite{BabaiBCIL96,CaiLZ00}.
 
 For the case of (in general) non-commutative groups, the knapsack problem has been studied in \cite{DudkinTreyer2018,FGLZ20,FrenkelNU15,GanardiKLZ18,KoenigLohreyZetzsche2015a,Loh19hyp,LohreyZ18,MyNiUs14}.
 In these papers, group elements are usually represented by finite words over the generators (although in \cite{LohreyZ18} a more succinct
 representation by so-called straight-line programs is studied as well). Note that for the group $\Z$ this corresponds to a unary representation
 of integers. Hyperbolic groups, which are of fundamental importance in the area of geometric group theory, are an important class of groups
 where knapsack can be decided in polynomial time (even in $\mathsf{LogCFL}$, i.e., the closure of the context-free languages under logspace
 reductions). This result can be extended to the class of all groups
 that can be built from hyperbolic groups by the operations of (i) direct products with $\Z$ and (ii) free products~\cite{LohreyZ18}. On the other hand,
 for many groups the knapsack problem is $\NP$-complete. Examples are certain right-angled Artin groups (like the direct product of two free groups
 of rank two~\cite{LohreyZ18}), wreath products (e.g. the wreath product $\Z \wr \Z$~\cite{GanardiKLZ18}) and 
 free solvable groups~\cite{FGLZ20}.  For wreath products $G \wr \Z$, where $G$ is finite non-solvable or free of rank at least two, 
 the knapsack problem is complete for $\Sigma^p_2$ (the second existential level of the polynomial time hierarchy) \cite{FGLZ20}.
 Finally, for finitely generated nilpotent groups, the knapsack problem is in general undecidable \cite{GanardiKLZ18,MishchenkoT17}.
 
Our second main result is that for the Baumslag-Solitar groups $\BS(1,q)$ with $q \geq 2$ the knapsack problem is $\NP$-complete
 (Theorem~\ref{thm-main-KP}). This extends a result from
 \cite{DudkinTreyer2018}, where decidability (without any complexity bound) was shown for a restriction of the knapsack problem for 
 $\BS(1,q)$. In this restriction, all group elements $g_i$ must be represented by words where the exponent sum of all occurrences of $t$ is zero
 (here we refer to the presentation $\langle a,t \mid t a t^{-1} = a^q \rangle$ of $ \BS(1,q)$). Showing $\NP$-hardness of 
 the knapsack problem for $\BS(1,q)$ is easy (based on the result that knapsack for $\Z$ with binary encoded integers is $\NP$-hard).
 For membership in $\NP$ we use a recent result of Gu\'{e}pin, Haase, and Worrell~\cite{GuepinHaaseWorrell2019} according to which the existential fragment of 
 B\"uchi arithmetic (an extension of Presburger arithmetic) belongs to $\NP$. The $\NP$-membership of the knapsack problem for $\BS(1,q)$
 is a bit of a surprise, since one can show that minimal solutions of knapsack equations over $\BS(1,q)$ can be of size doubly exponential in the 
 length of the equation, see Theorem~\ref{thm-2exp}. This rules out a simple guess-and-verify strategy.
 
\subsection{Solvability of systems of exponent equations.} 
In the final section of the paper we consider the following generalization of the knapsack problem for $\BS(1,q)$: the input is a conjunction
\begin{equation} \label{eq-exp-system}
\bigwedge_{i=1}^n g_{i1}^{x_{i1}} g_{i2}^{x_{i2}} \cdots g_{id_i}^{x_{id_i}} = h_i, 
\end{equation}
where the $g_{ij}, h_i$ are elements of $\BS(1,q)$ and the $x_{ij}$ are variables taking values in $\N$. In contrast to 
the knapsack problem, we do not require these variables to be pairwise different (we also allow $x_{ij} = x_{ik}$).
We call \eqref{eq-exp-system} a system of exponent equations. Aside from being a natural generalization of the knapsack problem,
systems of exponent equations play a crucial role in a characterization of decidability of the knapsack problem for wreath products~\cite{BergstrasserGZ21}: In order to understand for which wreath products $G\wr H$ the knapsack problem is decidable, we need to clarify for which groups $G$ one can decide solvability of systems of exponent equations. For example, the knapsack problem is decidable for $\Z\wr G$ if and only if solvability of systems of exponent equations is decidable for $G$ (this special case already follows from \cite[Proposition~3.1, Theorem~5.3]{GanardiKLZ18}).

For many groups, solvability of systems of exponent equations is
decidable. This holds in fact for all so-called knapsack semilinear groups, i.e., groups where the set of solutions of a knapsack equation
is an effectively computable semilinear set. Examples of knapsack semilinear groups are hyperbolic groups \cite{Loh19hyp}
and co-context-free groups \cite{KoenigLohreyZetzsche2015a}. Moreover, the class of knapsack semilinear groups is effectively 
closed under finite extensions \cite{FiLoZe21}, wreath products \cite{GanardiKLZ18}, graph products \cite{FiLoZe21}, and amalgamated products and HNN-extensions 
over finite groups \cite{FiLoZe21}. On the other hand, solvability of systems of exponent equations is undecidable for the discrete Heisenberg group 
\cite{KoenigLohreyZetzsche2015a}.

Our last main result states that solvability of systems of exponent equations is undecidable for every Baumslag-Solitar group $\BS(1,q)$ with $q \ge 2$ (Theorem~\ref{thm-main-EXP}).
We prove this result by a reduction from the existential theory of $(\N,+,(x,y) \mapsto x \cdot 2^y)$, which was shown to be undecidable
by B\"uchi and Senger \cite[Corollary~5]{BuchiS88}.
In contrast to this result, it has been shown recently that the Diophantine
theory (or, equivalently, solvability of systems of word equations with variables ranging over $\BS(1,q)$) is decidable for $\BS(1,q)$ \cite{KLM20}. 

A preliminary version of this paper appeared in \cite{LohreyZ20}.

\section{Preliminaries}
  
 For $a,b \in \Z$ we write $a \mid b$ ($a$ divided $b$) if $b = ka$ for some $k \in \Z$. 
 We denote with $[a,b]$ the interval $\{ z \in \Z \mid a \leq z \leq b\}$.
 For complex numbers $\alpha_1, \ldots, \alpha_k \in \C$ we denote with
 $\Z[\alpha_1, \ldots, \alpha_k]$ the subring of $\C$ obtained by adjoining 
 to the ring of integers $\Z$ the complex numbers  $\alpha_1, \ldots, \alpha_k$.
 
 The set of polynomials with variable $x$ and coefficients from $\Z$ is denoted with $\Z[x]$.
 Let $p(x)  = a_n x^n + a_{n-1} x^{n-1} + \cdots + a_1 x + a_0 \in \Z[x]$ with $a_n \neq 0$.
 Then we define  $\deg(p) = n$ (the {\em degree} of $p(x)$) and
 $\height(p) = \max\{ |a_0|, \ldots, |a_n|\}$. 
 The {\em dense representation} of the above polynomial is the tuple
 $(a_0, a_1, \ldots, a_n)$, where every $a_i$ is given in binary encoding.
 We define the {\em sparse representation} of a polynomial 
$p(x) = a_0 x^{e_0} + a_1 x^{e_1} + \cdots + a_n x^{e_n}$ 
with $a_i \in \Z \setminus \{0\}$ for all $0 \in [0,n]$ and $0 \leq e_0 < e_1 <\cdots < e_n$ 
as the list $(a_0, e_0, a_1, e_1,\ldots, a_n, e_n)$ where all numbers
in this list are written in binary representation.
 
 A \emph{Laurent polynomial} is a polynomial that may also contain
 powers $x^k$ with $k < 0$. Formally, a Laurent polynomial over $\Z$ is an expression
 $p(x) = \sum_{i \in \Z} a_i x^i$ with $a_i \in \Z$  such that only finitely many $a_i$
 are non-zero. With $\Z[x,x^{-1}]$ we denote the set of all Laurent 
 polynomials over $\Z$; it is a ring with the natural addition and multiplication operations. 
 If $p(x) = \sum_{i=k}^l a_i x^i$ with $k,l \in \Z$, $k \leq l$ and $a_k \neq 0 \neq a_l$ then we
 define the {\em dense unary} (resp., {\em dense binary}) {\em representation} of the Laurent polynomial $P(x)$ 
 as the list of unary (resp., binary) encoded integers $a_k, a_{k+1}, \ldots, a_l$ together with $k$ in unary encoding.
 
For a complex number $\alpha \in \C \setminus \{0\}$ we have a natural homomorphism from 
 $\Z[x,x^{-1}]$ to $\Z[\alpha,\alpha^{-1}]$ obtained by evaluating a Laurent polynomial
 at $x = \alpha$. Clearly, for an integer $q \in \Z\setminus \{0\}$ we have $\Z[q,q^{-1}] = \Z[1/q]$.
If $q \geq 2$, this is the set of all rational numbers that have 
finite expansion in base $q$, i.e., the set of all numbers
$\sum_{a \leq i \leq b} r_i q^i$ with $r_i \in [0,q-1]$ and $a,b \in \Z$.
If $u = \sum_{-k \leq i \leq \ell} r_i q^i \neq 0$ with $k, \ell \geq 0$ and
$\ell+k$ minimal, we define 
 $\|u\|_q = \ell+k+1$.  Under the assumption that $q$ is a constant
 (which will be always the case in this paper),  $\|u\|_q$ is the number of digits in the $q$-ary
 representation of $u$.

\subsection{Circuit complexity}  \label{sec-complexity}
 We assume basic knowledge in complexity theory, in particular with the complexity class $\NP$; see \cite{AroraBarak09} for details.
 We deal with the circuit complexity class $\TC^0$. It 
 contains all languages $L \subseteq \{0,1\}^*$ that can be solved by a family of threshold circuits of 
polynomial size and constant depth. More formally: we have a family $\mathcal{C} = (C_n)_{n \ge 0}$
of boolean circuits $C_n$ with the following properties:
\begin{itemize}
\item $C_n$ has exactly $n$ input gates $x_1, \ldots, x_n$ with fan-in zero (the fan-in of a gate is the number of incoming wires).
\item All other gates are either not-gates (with fan-in one), and-gates (with arbitrary fan-in), or majority-gates (with arbitrary fan-in). A majority gate of fan-in $k$ evaluates
to $1$ if and only if at least $k/2$ many input wires carry the truth value $1$.
\item Every $C_n$ has a distinguished output gate.
\item There is a constant $d$ such that the depth of very circuit $C_n$ is bounded by $d$, where the depth of a circuit is the length of a longest path from
an input gate to the output gate.
\item There is a polynomial $p(n)$ such that $C_n$ has at most $p(n)$ many gates.
\item For every word $w = a_1 a_2 \cdots a_n$ with $a_i \in \{0,1\}$, we have $w \in L$ if and only if the output gate of the circuit 
$C_n$ evaluates to $1$ when every input gate $x_i$ is set to $a_i$.
\end{itemize}
We can lift this definition to languages over an arbitrary alphabet $\Sigma$ by fixing a binary encoding of the symbols from $\Sigma$.
We always assume such encodings implicitly.
To compute a function $f \colon \{0,1\}^* \to \{0,1\}^*$ by a circuit family, we encode $f$ by the language
$L_f = \{1^i 0 w \mid w \in \{0,1\}^*, \text{the $i$-th bit of $f(w)$ is 1} \}$.

In this paper, we only deal with the $\mathsf{DLOGTIME}$-uniform version of $\TC^0$. In this variant, $\TC^0$ is contained in deterministic logspace and hence
in polynomial time.
We do not give the quite technical definition of $\mathsf{DLOGTIME}$-uniformity; see see \cite{Vol99} for details. In fact, all we need about $\TC^0$ is 
the fact  that the following problems
belong $\mathsf{DLOGTIME}$-uniform $\TC^0$: \label{sec-ct}
\begin{enumerate}
\item iterated addition/multiplication (i.e., addition/multiplication of an arbitrary number) of binary encoded numbers/polynomials that are given in dense representation  \cite{Eberly89,HesseAB02},
\item division with remainder of two binary encoded numbers/polynomials that are given in dense representation \cite{Eberly89,HesseAB02},
\item computing the number $|w|_a$ of occurrences of a letter $a$ in a word $w$,
\item computing an image $h(w)$ where $h : \Sigma^* \to \Gamma^*$ is a homomorphism of free monoids \cite{LaMc98}.
\end{enumerate}
The results on binary numbers hold for any basis, since one can transform between binary representation and 
$q$-ary representation; this is a consequence of the first two points.

In the rest of the paper, when we speak about $\TC^0$, we always refer to $\mathsf{DLOGTIME}$-uniform $\TC^0$.

\subsection{Algebraic numbers}

An algebraic number is a complex number which is the root of a polynomial from $\Z[x]$.
For every algebraic number $\alpha$ there is a unique polynomial $p(x) \in \Z[x]$
with $p(\alpha) = 0$ and such that $p(x)$ has minimal degree among all such polynomials
and the coefficients of $p(x)$ have no common divisor $>1$.
This polynomial is called the {\em minimal polynomial} of $\alpha$. If $p(x)$ is the minimal
polynomial of $\alpha$, then we define $\deg(\alpha) = \deg(p)$ and
$\height(\alpha) = \height(p)$.

{\em Sparse polynomial root testing} is the following 
decision problem:
\begin{description}
\item[Input] A polynomial $P(x) \in \Z[x]$ given in sparse representation and an algebraic number $\alpha \in \C$ given by its minimal polynomial in dense representation.
\item[Question] Is $P(\alpha) = 0$?
\end{description}
Note that we do not specify $\alpha$ uniquely:  if $p(x)$ is the minimal polynomial then by writing down only $p(x)$, we cannot distinguish
$\alpha$ from its conjugates. On the other hand, for sparse polynomial root testing there is no reason to make this distinction, because
$P(\alpha) = 0$ if and only if $P(x)$ is a multiple of $p(x)$.

\begin{theorem} \label{thm-root-testing}
Sparse polynomial root testing is in $\TC^0$.\footnote{We are not aware of a $\TC^0$-algorithm for testing whether a given polynomial is irreducible. Hence, in the statement of this theorem we consider sparse polynomial root testing as a promise problem. On the other hand, in our later application we will
only deal with a fixed algebraic number $\alpha$ in which case the minimal polynomial of $\alpha$ can be hard-wired in the algorithm.}

\end{theorem}

\begin{proof}
Let $p_\alpha(x)$ be the minimal polynomial of $\alpha$.
In \cite{Lenstra} it was shown that the problem belongs to polynomial time using the following gap theorem:
Let $P(x) = P_0(x) + x^s P_1(x) \in \Z[x]$ be a polynomial with $k+1$ monomials
and $u = \deg(P_0)$ and let $d \geq 1$ be an integer
such that 
\begin{equation} \label{eq-s-u}
s-u > \frac{\ln k + \ln\height(P)}{c(d)}
\end{equation}
where 
\[
c(1) = \ln 2 \quad\text{and}\quad c(d) = \frac{2}{d \cdot (\ln(3d))^3} \text{ for } d \geq 2.
\]
If $\alpha$ is an algebraic number of degree at most $d$ which is not a root of unity then
 $P(\alpha) = 0$ if and only if $P_0(\alpha)=0$ and $P_1(\alpha) = 0$. Note that the number on the right-hand side of 
 \eqref{eq-s-u} is polynomial in the input length if $P$ is given in spare representation and the minimal polynomial of $\alpha$ is given
 in dense representation.
This allows to split in $\TC^0$ the input polynomial $P(x)$ into several polynomials $p_0(x), \ldots, p_k(x)$  such that
$P(\alpha) = 0$ if and only if $p_i(\alpha)=0$ for all $1 \leq i \leq k$. Moreover, all $p_i$ are computed
in dense representation. Finally, we check for every $i$ whether $p_\alpha(x)$ divides $p_i(x)$.

It remains to consider the case where $\alpha$ is a root of unity. The case $\alpha = \pm 1$ is clear since
iterated addition is in $\TC^0$. Otherwise $\alpha$ is an $m^{\text{th}}$ primitive root of unity
for some $m > 2$ and the degree of $p_\alpha(x)$ is $d = \varphi(m)$, where $\varphi$
is Euler's phi-function.  We have $m \leq 3 \varphi(m)^{3/2} = 3 d^{3/2}$ \cite{BradfordD88}.
Hence, given $p_\alpha(x)$ of degree $d$ we simply test in parallel for every $d+1 \leq e \leq 3 d^{3/2}$ whether
$p_\alpha(x)$ divides $x^e-1$. Once we found such an $e$ we can replace in the polynomial $P(x)$ every
binary encoded monomial $x^n$ by $x^{n \bmod e}$. In this way we can compute a polynomial $\tilde P(x)$ 
in dense representation such that $\tilde P(\alpha) = 0$ if and only if  $P(\alpha) = 0$. Finally, we check
whether $p_\alpha(x)$ divides $\tilde P(x)$.
\end{proof}

\subsection{Groups}

 We assume that the reader is familiar with the basics of group theory.
Let $G$ be a group. We always write $1$ for the group identity element.
We say that $G$ is \emph{finitely generated (f.g.)} if there is a finite subset $\Sigma \subseteq G$
such that every element of $G$ can be written as a product of elements from $\Sigma$; such a $\Sigma$ is called
a \emph{(finite) generating set} for $G$. We always assume that $a \in \Sigma$ implies $a^{-1} \in \Sigma$;
such a generating set is also called \emph{symmetric}.
We write $G = \langle \Sigma \rangle$ if $\Sigma$ is a symmetric generating set for $G$. In this case, we have a canonical
surjective morphism $h : \Sigma^* \to G$ that maps a word over $\Sigma$ to its product in $G$ (the so called
{\em evaluation morphism}). If $h(w)=1$ we also say that $w=1$ \emph{in $G$}.
On $\Sigma^*$ we can define a natural involution $\cdot^{-1}$ by $(a_1 a_2 \cdots a_n)^{-1} = a_n^{-1} \cdots a_2^{-1} a_1^{-1}$
for $a_1, a_2, \ldots, a_n \in \Sigma$.

\subsubsection{Matrix groups}
 
For a complex number $\alpha \in \C \setminus \{0\}$ let $T(\alpha)$ be the subgroup of 
 $\GL(2,\C)$ consisting of the upper triangular matrices
\begin{equation} \label{matrix2}
\begin{pmatrix} \alpha^k & u \\ 0 & 1 \end{pmatrix}
\end{equation}
with $k\in\Z$ and $u\in\Z[\alpha,\alpha^{-1}]$.  This means we have the multiplication
\begin{equation} 
\begin{pmatrix} \alpha^k & u \\ 0 & 1 \end{pmatrix}\begin{pmatrix} \alpha^\ell & v \\ 0 & 1 \end{pmatrix}=\begin{pmatrix} \alpha^{k+\ell} & u+ \alpha^k \cdot v \\ 0 & 1\end{pmatrix}. \label{multiplication}
\end{equation}
This group can be also written as the semi-direct product $\Z[\alpha,\alpha^{-1}] \rtimes \Z$, where $\Z$ acts on $\Z[\alpha,\alpha^{-1}]$ by multiplication with $\alpha$.
The groups $T(\alpha)$ are also studied in \cite{Guyot12,Guyot18}.

We encode the matrix \eqref{matrix2} by the pair $(k,p)$, where $k$ is given in unary encoding and $p$ is a Laurent polynomial with $u = p(\alpha)$
 that is given in dense unary representation.
The group $T(\alpha)$ is generated by the two matrices
\begin{equation} \label{gl-A-T}
  a = \begin{pmatrix} 1 & 1 \\ 0 & 1 \end{pmatrix}  \quad\text{and}\quad  t = \begin{pmatrix} \alpha & 0 \\ 0 & 1 \end{pmatrix}
\end{equation}
and their inverses. We denote with $h \colon \{a,a^{-1}, t, t^{-1}\}^* \to T(\alpha)$ the canonical evaluation morphism.
Hence, $h(w)$ is the identity matrix if and only if $w = 1$ in $T(\alpha)$.

We now have two encodings of elements from $T(\alpha)$: as pairs $(k,p)$ describing a matrix  \eqref{matrix2} and as words over 
the alphabet $\{a,a^{-1}, t, t^{-1}\}$. By the the following lemma, we can switch in $\TC^0$ between these encodings.

\begin{lemma} \label{lemma-TC0-conversion}
Given a word $w \in \{a,a^{-1}, t, t^{-1}\}^*$ we can compute in $\TC^0$ the matrix $h(w)$ encoded as a pair $(k,p)$ as above.
Vice versa, given a matrix $A \in T(\alpha)$ in the above encoding, we can compute in $\TC^0$ a word $w \in h^{-1}(A)$.
\end{lemma}

\begin{proof}
First consider a word $w \in \{a,a^{-1}, t, t^{-1}\}^*$ and let $h(w)$ be the matrix in \eqref{matrix2}.
Then $k = |w|_t - |w|_{t^{-1}}$, which can be computed in $\TC^0$.
It remains to compute a Laurent polynomial $p(x)$ in dense unary representation such that $u = p(\alpha)$. 
Let $w_1 a^{\epsilon_1}, \ldots, w_l a^{\epsilon_l}$ be all prefixes of $w$
that end with $a$ or $a^{-1}$ ($\epsilon_1, \ldots, \epsilon_l \in \{-1,1\}$).
Let $k_i = |w_i|_t - |w_i|_{t^{-1}}$, which can be computed in $\TC^0$ in unary notation.
Then, $u = p(\alpha)$ with $p(x) = \sum_{i=1}^l \epsilon_i x^{k_i}$. The dense unary representation of this polynomial can be 
easily computed in $\TC^0$.

The inverse transformation is straightforward: take the matrix \eqref{matrix2}, where $k$ is given in unary encoding and 
$u = p(x)$ for a Laurent polynomial $p(x)$ in dense unary representation.
A matrix of the form $\gsmat{\alpha^z}{1}$
(for a unary encoded $z$) can be produced by the word $t^{z} a t^{-z}$. By concatenating such words (which is possible in $\TC^0$ by
point 4 from page~\pageref{sec-ct}), one can produce
from a given Laurent polynomial $p(x)$ in dense unary representation
 a word for the matrix $\gsmat{p(\alpha)}{1}$. Finally, one has to concatenate 
$t^k$ on the right in order to produce the matrix  \eqref{matrix2}.
\end{proof}

\subsubsection{Baumslag-Solitar groups}

For $p,q \in \Z \setminus \{0\}$,
 the \emph{Baumslag-Solitar group} $\BS(p,q)$ is defined as the finitely presented group
 $\BS(p,q) = \langle a,t \mid t a^p t^{-1} = a^q \rangle$.
 We can w.l.o.g.~assume that $q \geq 1$.
 Of particular interest are the Baumslag-Solitar groups $\BS(1,q)$ for $q \geq 2$.
 They are solvable and linear.
It is well-known (see e.g. \cite[III.15.C]{Woe00}) that $\BS(1,q)$ is isomorphic to $T(q)$.
Moreover, the generator $a$ (resp., $t$) of $\BS(1,q)$ corresponds to the matrix $a$ (resp., $t$)
from \eqref{gl-A-T}.
From Lemma~\ref{lemma-TC0-conversion} we immediately get:

\begin{lemma} \label{lemma-TC0-conversion2}
Given a word $w \in \{a,a^{-1}, t, t^{-1}\}^*$ we can compute in $\TC^0$ the matrix $h(w)$ with matrix entries given in $q$-ary encoding.
Vice versa, given a matrix $A \in T(q)$ with $q$-ary encoded entries, we can compute in $\TC^0$ a word $w \in h^{-1}(A)$.
\end{lemma}
By the previous lemma, we can represent elements of $\BS(1,q)$ either as words over the alphabet $\{a,a^{-1}, t, t^{-1}\}$
or by matrices from $T(q)$ with $q$-ary encoded entries. 
For the matrix $A\in T(q)$ in \eqref{matrix2} (with $\alpha=q$) we define $\|A\| = |k| + \| u \|_q$. 
Hence $\|A\|$ is the length of the encoding of $A$.
  
  Another well known special case of the group $T(\alpha)$ is obtained when $\alpha$ is transcendental. In this case $T(\alpha)$ is 
  isomorphic to the wreath product $\Z \wr \Z$:  It is isomorphic to the group of all matrices
\begin{equation} \label{eq-matrix-ZwrZ}
\begin{pmatrix}
  x^k & P(x) \\
  0 & 1 
\end{pmatrix}
\end{equation}
where $k \in \Z$ and $P(x) \in \Z[x,x^{-1}]$  (see e.g.~\cite[Section~2.2]{MRUV10}).
In contrast to  $\BS(1,q)$ the group $\Z \wr \Z$ is not finitely presented~\cite{baumslag1961wreath}.
A well-known infinite presentation of $\Z \wr \Z$ is $\langle a,t \mid [a^{t^i}, a^{t^j}] = 1 \ (i,j \in \Z)\rangle$.

\subsubsection{Knapsack, exponent equations and the power word problem} \label{sec-knapsack}
 
 Let $G = \langle \Sigma \rangle$ be a f.g.~group.
Moreover, let $x_1, x_2, \ldots, x_d$ be pairwise distinct variables.
A \emph{knapsack expression} over $G$ is an expression of the form 
\begin{equation} \label{ks-expression}
E = v_0 u_1^{x_1} v_1 u_2^{x_2} v_2  \cdots u_d^{x_d} v_d
\end{equation}
with $d \geq 1$, words $v_0,\ldots, v_d \in \Sigma^*$ and non-empty words $u_1, \ldots, u_d \in \Sigma^*$.
 A tuple $(n_1, \ldots, n_d) \in\N^d$ is a  \emph{$G$-solution} of $E$ if $v_0 u_1^{n_1} v_1 u_2^{n_2} v_2  \cdots u_d^{n_d} v_d = 1$ in $G$.
 With $\Sol(G,E)$ we denote the set of all $G$-solutions of $E$.   
The \emph{size} of $E$ is defined as $|E| = |v_0|+\sum_{i=1}^d |u_i|+|v_i|$.
The {\em knapsack problem for $G$}, $\KP(G)$ for short, is the following 
decision problem:
\begin{description}
\item[Input] A knapsack expression $E$ over $G$.
\item[Question] Is $\Sol(G,E)$ non-empty?
\end{description}
It is easy to observe that the concrete choice of the generating set $\Sigma$ has no influence
on the decidability/complexity status of $\KP(G)$.
W.l.o.g.~we can restrict to knapsack expressions of the form $u_1^{x_1} u_2^{x_2}  \cdots u_d^{x_d} v$:
 for $E = v_0 u_1^{x_1} v_1 u_2^{x_2} v_2   \cdots u_d^{x_d} v_d$
and 
\[ E' = (v_0 u_1 v_0^{-1}) ^{x_1} (v_0 v_1 u_2 v_1^{-1} v_0^{-1})^{x_2} \cdots (v_0 \cdots v_{d-1} u_d v_{d-1}^{-1} \cdots v_0^{-1})^{x_d}
v_0 \cdots v_{d-1} v_d
 \]
we have $\Sol(G,E) = \Sol(G,E')$.

An {\em exponent expression} over $G = \langle\Sigma\rangle$ is a formal 
expression $E$ as in \eqref{ks-expression}, but in contrast to knapsack expressions, we allow $x_i = x_j$ for $i \neq j$. 
The set of solutions $\Sol(G,E)$ for the exponent expression $E$ can be defined analogously to knapsack expressions.
We define {\em solvability of systems of exponent equations over $G$}, $\ExpEq(G)$ for short,  as the following decision problem:
\begin{description}
\item[Input] A finite list of exponent expressions $E_1,\ldots,E_n$ over $G$.
\item[Question] Is $\bigcap_{i=1}^n \Sol(G,E_i)$ non-empty?
\end{description}
This problem has been studied for various groups in \cite{FGLZ20,GanardiKLZ18,Loh19hyp,LohreyZ18}.

A \emph{power word} (over $\Sigma$) is a tuple $(u_1,k_1,u_2,k_2,\ldots,u_d,k_d)$ where 
$u_1, \dots, u_d \in \Sigma^*$ are words over the group generators 
and $k_1, \dots, k_d\in \mathbb{Z}$ are integers that are given in binary notation. Such a power word represents the 
word $u_1^{k_1} u_2^{k_2}\cdots u_d^{k_d}$. Quite often, we will identify the power word $(u_1,k_1,u_2,k_2,\ldots,u_d,k_d)$
with the word $u_1^{k_1} u_2^{k_2}\cdots u_d^{k_d}$. 
The \emph{power word problem} for the f.g.~group $G$, $\PowWP(G)$ for short, is defined as follows:
\begin{description}
\item[Input] A power word $(u_1,k_1,u_2,k_2,\ldots,u_d,k_d)$.
\item[Question] Does $u_1^{k_1} u_2^{k_2}\cdots u_d^{k_d}=1$ hold in $G$?
\end{description}
Due to the binary encoded exponents, a power word can be seen as a succinct description 
of an ordinary word. The size of the above power word $w$ is 
$\sum_{i=1}^d |u_i| + \lceil\log_2 k_i\rceil$ which is the length of the binary encoding of $w$.

 \section{Power word problem for BS(1,q)}

In this section we prove our first main result:
\begin{theorem}  \label{thm-powerWP}
For every $\alpha \in \C \setminus \{0\}$,
$\PowWP(T(\alpha))$ belongs to $\TC^0$.
\end{theorem}

\begin{proof}
If $\alpha$ is transcendental, then $T(\alpha)$ is isomorphic to $\Z \wr \Z$ for which the power word problem belongs
to $\TC^0$ \cite{LoWe19}. For the rest of the proof we assume that $\alpha$ is algebraic.
We show that in this case, $\PowWP(T(\alpha))$ is $\TC^0$-reducible to sparse polynomial root testing, which belongs to $\TC^0$
by Theorem~\ref{thm-root-testing}.

Let us fix an algebraic number $\alpha \in \C\setminus \{0\}$.
Consider a power word of the form
\[
\begin{pmatrix} \alpha^{k_1} & u_1 \\ 0 & 1 \end{pmatrix}^{\!\! n_1} \cdot \begin{pmatrix} \alpha^{k_2} & u_2 \\ 0 & 1 \end{pmatrix}^{\!\! n_2} \cdots \begin{pmatrix} \alpha^{k_l} & u_l \\ 0 & 1 \end{pmatrix}^{\!\! n_l}
\]
Here, the $n_i$ are binary encoded integers  and every $u_i$ is of the form $u_i = p_i(\alpha)$ for a Laurent polynomial $p_i$ over $\Z$ that is given
in dense unary representation. Note that
\[
\gmat{0}{\alpha^{m}} \gmat{u}{\alpha^{n}} \gmat{0}{\alpha^{-m}} = \gmat{\alpha^m u }{\alpha^{n}} .
\]
By this we can assume that all $p_i$ are ordinary polynomials over $\Z$.
We have
\begin{eqnarray*}
\begin{pmatrix} \alpha^{k_i} & u_i \\ 0 & 1 \end{pmatrix}^{\!\! n} &=& 
\gmat{(1+\alpha^{k_i} + \alpha^{2k_i} + \cdots + \alpha^{(n-1)k_i}) \cdot u_i}{\alpha^{k_i n}} \\
&=&
\begin{cases}
\gmat{\frac{\alpha^{k_i n} - 1}{\alpha^{k_i} - 1} \cdot u_i}{\alpha^{k_i n}} & \text{ if $\alpha^{k_i}\neq 1$} \\
\gmat{n \cdot u_i}{1} & \text{ if $\alpha^{k_i}=1$}
\end{cases}
\end{eqnarray*}
Hence, our power word can be written as
\begin{eqnarray*}
&& \gmat{v_1}{\alpha^{k_1 \cdot n_1}} \cdot \gmat{v_2}{\alpha^{k_2 \cdot n_2}} \cdots \gmat{v_l}{\alpha^{k_l \cdot n_l}} \\
&=&
\gmat{v_1 + \alpha^{k_1 n_1} v_2 + \cdots + \alpha^{k_1 n_1 + \cdots +k_{l-1} n_{l-1}} v_l}{\alpha^{k_1 n_1 + \cdots + k_l n_l}}
\end{eqnarray*}
Here, $v_i$ is $n_i \cdot p_i(\alpha)$ (if $\alpha^{k_i}=1$)
or $\frac{\alpha^{k_i n_i} - 1}{\alpha^{k_i} - 1} \cdot p_i(\alpha)$  (if $\alpha^{k_i}\neq 1$).

We have to check whether   
\begin{eqnarray} 
\alpha^{k_1 n_1 + \cdots + k_l n_l} & =& 1 \label{eq=-1} \\
\label{eq=0}
v_1 + \alpha^{k_1 n_1} v_2 + \cdots + \alpha^{k_1 n_1 + \cdots +k_{l-1} n_{l-1}} v_l &=& 0.
\end{eqnarray} 
Equality \eqref{eq=-1} can be easily checked in $\TC^0$: we compute in $\TC^0$ the binary encoding of $s = k_1 n_1 + \cdots + k_l n_l$.
If $\alpha$ is not a root unity then we check 
whether $s=0$. On the other hand, iff $\alpha$ is a primitive, say $d^{\text{th}}$, root of unity, then
we check whether $d$ divides $s$.

The verification of \eqref{eq=0} can be reduced to sparse polynomial root testing as follows.
First, we compute all binary encoded numbers $s_i = k_1 n_1 + \cdots + k_{i-1} n_{i-1}$ for
$i \in [1,l+1]$. By multiplying \eqref{eq=0} with a power $\alpha^m$ for $m \geq 0$ sufficiently
large, we can assume that all $s_i$ are non-negative. We have to check whether
\begin{equation} \label{eq=1}
\sum_{i=1}^l  \alpha^{s_i} v_i = 0.
\end{equation} 
Let $J = \{ i \in [1,l] \mid \alpha^{k_i} \neq 1\}$
and define the polynomial
\[
q(x) = \prod_{i \in J} (x^{k_i} - 1) .
\]
Note that $q(\alpha) \neq 0$.
We can compute in $\TC^0$ the dense representation of $q(x)$ (recall from Section~\ref{sec-complexity} that iterated multiplication of densely represented polynomials
is in $\TC^0$). Then, we compute for all $i \in [1,l]$ the sparse representation of the polynomial
\[
q_i(x) := \begin{cases}
n_i \cdot q(x) \cdot p_i(x) \cdot x^{s_i} & \text{ if $\alpha^{k_i}= 1$} \\
(x^{k_i n_i}-1) \cdot \prod_{j \in J \setminus \{i\}} (x^{k_j}-1) \cdot p_i(x) \cdot x^{s_i} & \text{ if $\alpha^{k_i}\neq 1$}.
\end{cases}
\]
This is possible in $\TC^0$. For instance, in the second case ($\alpha^{k_i}\neq 1$), we first compute in $\TC^0$ the 
 the dense representation of $\prod_{j \in J \setminus \{i\}} (x^{k_j}-1) \cdot p_i(x)$ (this is  iterated multiplication of densely represented polynomials).
 This dense representation can be easily multiplied in $\TC^0$ with the sparse representation of $(x^{k_i n_i}-1) \cdot x^{s_i} = x^{s_{i+1}}- x^{s_i}$, which yields the sparse
 representation of $q_i(x)$.

Finally we compute in $\TC^0$ the sparse representation of the polynomial.
\[ 
Q(x) = \sum_{i=1}^l q_i(x).
\]
We obtain 
\[
Q(\alpha) =  q(\alpha) \cdot \sum_{i=1}^l  \alpha^{s_i} v_i  .
\]
Since $q(\alpha) \neq 0$,  \eqref{eq=1} is equivalent to $Q(\alpha)  = 0$.
This concludes our reduction to sparse polynomial root testing.
\end{proof}

\section{Knapsack for BS(1,q)}

Whether the knapsack problem is decidable for $\BS(1,q)$ was left open
in \cite{DudkinTreyer2018}. Our second main result gives a positive answer 
and also settles the computational complexity:

\begin{theorem} \label{thm-main-KP}
For every $q \geq 2$,
$\KP(\BS(1,q))$ is $\NP$-complete.
\end{theorem}
Let us first remark that $\BS(1,q)$ is unusual in terms of its
knapsack solution sets.  In almost all groups where knapsack is known to
be decidable, knapsack equations have semilinear solution
sets~\cite{FGLZ20,FiLoZe21,GanardiKLZ18,KoenigLohreyZetzsche2015a,Loh19hyp,LohreyZ18}.
After the discrete Heisenberg group~\cite{KoenigLohreyZetzsche2015a},
the groups $\BS(1,q)$ are only the second known example where this is
not the case: the knapsack equation
$t^{-x_1}a^{x_2}t^{x_3}=a$ has the non-semilinear solution set
$\{(k,q^k,k) \mid k\in\N\}$.

Another unusual aspect is that knapsack is in $\NP$ although there are
knapsack equations over $\BS(1,2)$ whose solutions are all at least
doubly exponential in the size of the equation:
\begin{theorem}  \label{thm-2exp}
There is a family $E_k = E_k(x,y,z)$, $k \geq 1$, of solvable knapsack expressions over $\BS(1,2)$ such that
$|E_k| = \Theta(k)$ and $z \geq (2^{2\cdot 3^{k-1}}-1)/3^k-1$ for every solution of $E_k = 1$.
\end{theorem}

\begin{proof}
It is a well-known fact in elementary
number theory that for every $k \ge 1$, $2$ is a primitive root modulo~$3^k$, i.e., $2$ generates the group
$(\Z/3^k\Z)^*$ (the group of units of $\Z/3^k\Z$). See, for example, Theorem~3.6 and the remarks before
Theorem~3.8 in~\cite{Nathanson2000}. Consider the knapsack equation
\begin{equation} 
\gmat{0}{2}^x \gmat{1}{1} \gmat{0}{2^{-1}}^{y} \gmat{-3^{k}}{1}^z = \gmat{3^k+1}{1} \label{double-exponential-example}\end{equation}
in $\BS(1,2)$. 
In the top-left entry, it implies $2^{x} 2^{-y} = 1$.
Therefore, we must have
$x=y$ in every solution. In this case, the left-hand side of \cref{double-exponential-example} is
\[ \gmat{0}{2^{x}}\gmat{1}{1}\gmat{0}{2^{-x}}\gmat{-z3^{k}}{1} 
=\gmat{2^{x}-z\cdot 3^k}{1}.\]
Therefore, \cref{double-exponential-example} is equivalent to $x=y$
and $2^{x}-z\cdot 3^k=3^k+1$. Since some non-zero power of $2$ is congruent to $1$
modulo $3^k$, \cref{double-exponential-example} has a
solution. Moreover, any solution must satisfy $2^{x}\equiv 1\pmod{3^k}$. Since $2$ is a primitive root modulo $3^k$, $x$ must be a multiple of $|(\Z/3^k\Z)^*| = \varphi(3^k) = 2 \cdot 3^{k-1}$
(here, $\varphi$ is Euler's phi-function). Moreover, $x$ must be non-zero, because $1-z\cdot 3^k=3^k+1$ is not possible for $z\in\N$.
We obtain $x\ge 2\cdot 3^{k-1}$. Since $2^{x}-z\cdot 3^k=3^k+1$, this yields $z=(2^{x}-3^k-1)/3^k\ge (2^{2\cdot 3^{k-1}}-1)/3^k-1$. 
\end{proof}
\begin{remark}
  Subject to Artin's conjecture on primitive roots~\cite{Hooley1967}, a similar
  doubly-exponential lower bound results for every $\BS(1,q)$ where $q\ge 2$ is not
  a perfect square. Moreover, \Cref{thm-2exp} holds even if the
  variables $x,y,z$ range over $\Z$. For this, one replaces $3^k+1$
  with the inverse of $2$ in $(\Z/3^k\Z)^*$ in \cref{double-exponential-example}.
\end{remark}
Theorem~\ref{thm-2exp} rules out a simple guess-and-verify strategy to show Theorem~\ref{thm-main-KP}.
If one has an exponential upper bound (in terms of input length)
on the size of a smallest solution of a knapsack equation, then one can guess the binary
representation of a solution and verify, using the power word problem, whether the guess is indeed a solution.
The second step (verification of a solution using the power word problem) 
would work for $\BS(1,q)$ in polynomial time due to Theorem~\ref{thm-powerWP},
but the first step (guessing a binary encoded candidate for a solution) does not work for $\BS(1,2)$ due to
Theorem~\ref{thm-2exp}.

Our main tool for the proof of Theorem~\ref{thm-main-KP} is a recent result from 
 \cite{GuepinHaaseWorrell2019} concerning the existential fragment of B\"{u}chi arithmetic.

\subsection{B\"{u}chi arithmetic}
\emph{B\"{u}chi arithmetic}~\cite{buchi1960weak} 
 is the first-order theory of the structure $(\Z,+,\ge,0,V_q)$.
Here, $V_q$ is the function that maps $n\in\Z$ to the largest power of
$q$ that divides $n$. It is well-known that B\"{u}chi arithmetic is
decidable (this was first claimed in \cite{buchi1960weak}; a correct
proof was given in~\cite{Bruyere1985}). We will rely on the following
recent result of Gu\'{e}pin, Haase, and Worrell \cite{GuepinHaaseWorrell2019}:

\begin{theorem}[c.f.~\cite{GuepinHaaseWorrell2019}] \label{thm-existential-Buechi}
The existential fragment of B\"{u}chi arithmetic belongs to $\NP$.
\footnote{The paper~\cite{GuepinHaaseWorrell2019} shows an $\NP$
  upper bound for the structure $(\N,+,0,V_q)$, but an existential 
  sentence over the structure $(\Z,+,\ge,0,V_q)$ easily translates into one over
  $(\N,+,0,V_q)$.} 
\end{theorem}
We will also make use of the following simple lemma:

\begin{lemma} \label{lemma-mult}
Given the $q$-ary representation of a number $r \in \Z[1/q]$ we can construct in polynomial time
an existential Presburger formula over $(\Z,+)$ of size $\mathcal{O}(\| r \|_q)$ which expresses $y = r \cdot x$ for $x,y \in \Z$. 
\end{lemma}

\begin{proof}
Let $r = \sum_{-k \leq i \leq \ell} a_i q^i$ with $k,\ell \geq 0$ and $0 \leq a_i < q$ for $-k \leq i \leq \ell$.
We have $y = r x$ if and only if $q^k y = r' x$ for $r' = \sum_{i=0}^{k+\ell} a_{i-k} q^i \in \Z$.
Using iterated multiplication with the constant $q$ (which can be replaced by addition) we can easily define 
from $x$ and $y$ the integers $q^k y$ and $r' x$ by Presburger formulas of size $\mathcal{O}(k+\ell) = \mathcal{O}(\| r \|_q)$.
\end{proof}

\subsection{Proof of Theorem~\ref{thm-main-KP}}

We start with the lower bound.
The \emph{multisubset sum problem} asks for
integers $a_1,\ldots,a_d,b\in\Z$ given in binary, whether there exist
natural numbers $x_1,\ldots,x_d\ge 0$ with
$x_1a_1+\cdots+x_da_d=b$. It is known to be
$\NP$-complete~\cite{Haa11}. Since the knapsack equation
\[ \gmat{a_1}{1}^{x_1}\cdots \gmat{a_d}{1}^{x_d}=\gmat{b}{1} \]
is equivalent to $x_1a_1+\cdots+x_da_d=b$, we obtain $\NP$-hardness of knapsack over $\BS(1,q)$.
Note that computing the $q$-ary representation of $a_i$ from the binary representation is possible in logspace
(even in $\TC^0$).

For the upper bound we reduce $\KP(\BS(1,q))$ to the existential fragment of B\"{u}chi arithmetic, which belongs to $\NP$
by Theorem~\ref{thm-existential-Buechi}. We proceed in three steps.

\medskip

\paragraph{\bf Step 1:}  Expressing $M_g$ and $M_g^*$ using $S_\ell$.
We first express a particular set of binary relations using
existential first-order formulas over
$(\Z, +, \ge, 0, V_q, (\shift{q}{\ell})_{\ell\in\Z})$. Here, for
$\ell\in\Z$, $\shift{q}{\ell}$ is the binary predicate with
\[ x\,\shift{q}{\ell}\,y \iff \exists r\in\N \; \exists s\in\N\colon x=q^{r} \wedge y=q^{r+\ell\cdot s}. \]
Let $T_\Z(q)$ denote the subset of matrices in $T(q)$ that have entries in
$\Z$. We represent the matrix $\gsmat{n}{m} \in T_\Z(q)$ by the pair $(m,n)\in\Z\times\Z$
(note that we must have $m \in \N$).
Observe that we can define in the structure $(\Z, +, \ge, 0, V_q, (\shift{q}{\ell})_{\ell\in\Z})$ the set of pairs $(m, n)\in\Z$ such that
$\gsmat{n}{m}\in T_\Z(q)$, because this is equivalent to $m$ being a power
of $q$, which is expressed by  $1\,\shift{q}{1}\, m$.

A key trick is to express solvability of a knapsack equation
$g_1^{x_1}\cdots g_d^{x_d}=g$ without introducing variables in the
logic for $x_1,\ldots,x_d$. Instead, we employ the following binary relations
$M_g$ and $M_g^*$ on $T_\Z(q)$, which allow us to express existence of
powers implicitly. For $g\in T(q)$ and $x,y\in T_\Z(q)$, we have:
\begin{itemize}
\item $x\,M_g\, y \iff y=xg$,
\item $x\,M^*_g\,y \iff \exists s\in\N\colon y=xg^s$.
\end{itemize}
We construct existential formulas of size polynomial in $\|g\|$ over  the structure
$(\Z, +, \ge, 0, V_q, (\shift{q}{\ell})_{\ell\in\Z})$, which define the relations $M_g$ and $M_g^*$.
For the further consideration let
\[ g=\gmat{v}{q^\ell}.\]
Note that the relation $M_g$ is easily expressible because we
can express multiplication with $q^\ell$ and $v$ by existential Presburger formulas of length $\|g\|$, see Lemma~\ref{lemma-mult}.

We now focus
on the relations $M^*_g$ and express
\begin{equation} \label{eq-M*}
\gmat{u}{q^k}M^*_g \gmat{w}{q^m}
\end{equation}
Observe that for $\ell\ne 0$, we have
\begin{align*}
  \gmat{u}{q^k}\gmat{v}{q^\ell}^s &= \gmat{u}{q^k}\gmat{v+q^\ell v + \cdots + q^{(s-1)\ell} v}{q^{\ell s}} \\
                                  &=\gmat{u}{q^k}\gmat{v\frac{q^{\ell s}-1}{q^\ell-1}}{q^{\ell s}} = \gmat{u+v\frac{q^{k+\ell s} - q^k }{q^{\ell} -1 }}{q^{k+\ell s}}.
\end{align*}
Therefore, \cref{eq-M*} is equivalent to
\[ \exists x\in\Z \; \exists s\in\N\colon q^m= q^{k+\ell s} \wedge w=u+vx \wedge (q^\ell-1)x=q^m-q^k.\]
Here, we can quantify $x$ over $\Z$, because 
$$\frac{q^{k+\ell s}-q^k}{q^\ell-1} = q^k+q^{k+\ell}  + \cdots + q^{k+(s-1)\ell}$$
must be an integer
($k$ and $k+\ell s = m$ are non-negative).
Note that since we can express multiplication with $v$ and $q^\ell$ by existential Presburger formulas of size 
$\mathcal{O}(\|g\|)$ (Lemma~\ref{lemma-mult}), we can also express $w=u+vx$
and $(q^\ell-1)x=q^m-q^k$ by existential Presburger formulas of size $\mathcal{O}(\|g\|)$. Finally, we can express
$\exists s\in\N\colon q^m=q^{k+\ell s}$ using $q^k\,\shift{q}{\ell}\,q^m$.

It remains to express \cref{eq-M*} in the case
$\ell=0$. Note that 
$$g^s=\gmat{sv}{1}$$ in this case. Therefore, \cref{eq-M*} is equivalent to 
\begin{enumerate}[(i)]
\item there exists $s\in\N$ with $w=u+q^k\cdot s\cdot v$ and
\item $q^m=q^k$.
\end{enumerate}

Note that condition (i) is equivalent to
$\exists t\in\N\colon V_q(t)\ge q^k\wedge w=u+v\cdot t$. This is because choosing
$t=q^k\cdot s$ yields (i).  By Lemma~\ref{lemma-mult}, $w=u+v\cdot t$ can be expressed
by an existential Presburger formula of size $\mathcal{O}(\|g\|)$.

\medskip

\paragraph{\bf Step 2:} Expressing $S_\ell$ using $V_q$.
In our second step, we show that the binary relations $M_g$ and
$M_g^*$ can be expressed using existential formulas over
$(\Z,+,\ge,0,V_q)$ of size $\poly(\|g\|)$. As shown above, for this it suffices to define
$S_\ell$ by an existential formula over $(\Z,+,\ge,0,V_q)$ of size $\poly(\ell)$
 (note that the relations $S_\ell$ occur only positively
in the formulas from Step~1). For $m\in\N$, let $P_m$ be the
predicate where $P_m(x)$ states that $x$ is a power of $m$. We first
claim that for each $\ell\ge 0$, we can express $P_{q^\ell}$ using an
existential formula of size polynomial in $\ell$ over
$(\Z,+,\ge,0,V_q)$. The case $\ell=0$ is clear. 
For the case $\ell \ge 1$ we use the following observation from the proof of Proposition~7.1 in
\cite{bruyere1994logic}. Note that
$P_q(x)$ is just $V_q(x)=x$. 
\begin{fact} For all $\ell \ge 1$, $P_{q^\ell}(x)$ if and only if $P_q(x)$ and $q^\ell-1$ divides $x-1$.
\end{fact}
\begin{proof}
  If $x$ is a power of $q^{\ell}$, then $x=q^{\ell\cdot s}$ for some $s \geq 0$. So, $x$ 
  is a power of $q$. Moreover,
  \[ \frac{x-1}{q^{\ell}-1}=\frac{q^{\ell\cdot s}-1}{q^\ell-1}=\sum_{i=0}^{s-1}
  q^{i \ell}\] is an integer.

  Conversely, suppose $x$ is a power of $q$ and $q^{\ell}-1$ divides
  $x-1$. Write $x=q^{\ell\cdot s+r}$ with $0\le r<\ell$. Observe that
   \[x-1=q^{s\ell+r}-1=q^r(q^{s\ell}-1)+(q^r-1).\]
  Since $q^{\ell}-1$ divides $x-1$ as well as $q^{s\ell}-1$, we
  conclude that $q^{\ell}-1$ divides $q^r-1$. As $0\le r<\ell$,
  this is only possible with $r=0$. This shows the above fact.
\end{proof}
Using the predicates $P_{q^\ell}$, we can now express $S_\ell$.
Note that for $\ell\ge 0$, we have $x\,S_\ell \,y$ if and only if
\[ y\ge x \wedge \bigvee_{i=0}^{\ell-1} P_{q^\ell}(q^i x) \wedge P_{q^\ell}(q^i y).\]
Furthermore, for $\ell<0$, we have $x \,S_\ell \,y$ if and only if
$y \,S_{|\ell|} \,x$. Therefore, we can express each $S_\ell$ using
an existential formula of size polynomial in $\ell$ over
$(\Z,+,\ge,0,V_q)$. Hence, we can express $M_g$ and $M_g^*$ using
existential formulas of size $\poly(\|g\|)$ over $(\Z,+,\ge,0,V_q)$.

\medskip

\paragraph{\bf Step 3:} Expressing solvability of knapsack.
In the last step, we express solvability of a knapsack equation by an
existential first-order sentence over $(\Z,+,\ge,0,V_q)$, using the predicates
$M_g$ and $M_g^*$.  We claim that $g_1^{x_1}\cdots g_d^{x_d}=g$ has a
solution $(x_1,\ldots,x_d)\in\N^d$ if and only if there exist
$h_0,\ldots,h_{d}\in T_\Z(q)$ with
\begin{equation} 
h_0 M_{g_1}^*h_1 \ \wedge \ h_1M_{g_2}^* h_2 \ \wedge \cdots \wedge \ h_{d-1}M_{g_d}^* h_d \ \wedge \ h_0 M_{g} h_{d}. \label{knapsack-condition}\end{equation}
This can be stated by an existential sentence over $(\Z,+,\ge,0,V_q)$ of size polynomial in 
$\|g\|+\sum_{i=1}^d \|g_i\|$.

If such $h_0,\ldots,h_d$ exist, then for some $x_1,\ldots,x_d\in\N$,
we have $h_i=h_{i-1}g_i^{x_i}$ for all $i \in [1,d]$ and $h_d=h_0g$, which implies
$g_1^{x_1}\cdots g_d^{x_d}=g$. For the converse, we observe that for
each matrix $A\in T(q)$, there is some large enough $k\in\N$ such that
$\gsmat{0}{q^k}A \in T_\Z(q)$. Therefore, if
$g_1^{x_1}\cdots g_d^{x_d}=g$, then there is some large enough
$k\in\N$ so that for every $i \in [1,d]$, the matrix
$\gsmat{0}{q^k} g_1^{x_1}\cdots g_i^{x_i}$ has integer entries. With
this, we set $h_0=\gsmat{0}{q^k}$ and $h_i=h_{i-1}g_i^{x_i}$ for
$i \in [1,d]$.  Then we have $h_0,\ldots,h_d\in T_\Z(q)$ and
\cref{knapsack-condition} is satisfied.
\qed

\section{Systems of exponent equations over BS(1,q)}

Our algorithm for the knapsack problem in $\BS(1,q)$ cannot be extended
to solvability of systems of exponent equations (not even to solvability of a single exponent equation).
If we allow systems of exponent equations, we can show undecidability:

\begin{theorem} \label{thm-main-EXP}
For every $q \in \N$ with $q \geq 2$, $\ExpEq(\BS(1,q))$ is undecidable.
\end{theorem}

\begin{proof}
Consider the function $(x,y) \mapsto x \cdot 2^y$ on the natural numbers. B\"uchi and Senger \cite[Corollary~5]{BuchiS88} have shown that the existential
fragment of the first-order theory of $(\mathbb{N}, +, x \cdot 2^y)$ is undecidable. The proof generalizes to every function $(x,y) \mapsto x \cdot q^y$ 
for $q \in \N$, $q \geq 2$.  We reduce this fragment to $\ExpEq(\BS(1,q))$. For this it suffices to consider an existentially quantified conjunction of formulas of the 
following form: $x \cdot q^y = z$, $x+y = z$, and $x < y$ (the latter allow to express inequalities).
We replace each of these formulas by an equivalent exponent equation over $\BS(1,q)$.
For this we use the two generators $a$ and $t$ from \eqref{gl-A-T} (for $\alpha=q$).

The formula $x+y=z$ is clearly equivalent to $a^x a^y = a^z$, i.e., $a^x a^y a^{-z} = 1$. The formula  $x < y$ is equivalent 
$\exists z \in \N \colon a^x a^z a \, a^{-y}=1$. Finally, $x \cdot q^y = z$ is equivalent to $t^y a^x t^{-y} a^{-z}=1$.
\end{proof}

\section{Open problems}

Several open problems arise from our work:
\begin{itemize}
\item What is the complexity/decidability status of the power word/knapsack problem for Baumslag-Solitar groups
$\BS(p,q) = \langle a,t \mid t a^p t^{-1} = a^q \rangle$ for $p,q \geq 2$? Decidability of knapsack in case $\mathrm{gcd}(p,q)=1$
was shown in \cite{DudkinTreyer2018}, but the complexity as well as the decidability in case $\mathrm{gcd}(p,q)>1$ are open.
Since the word problem for $\BS(p,q)$ can be solved in logspace \cite{Weiss16}, one can easily show that
the power word problem for $\BS(p,q)$ belongs to $\mathsf{PSPACE}$.
By using techniques from \cite{LoWe19} one might try to find a logspace reduction from the power word problem for $\BS(p,q)$
to the word problem for $\BS(p,q)$ (the same was done for a free group in \cite{LoWe19}); this would show that the power word problem for $\BS(p,q)$
can be solved in logspace.
\item Baumslag-Solitar groups $\BS(1,q)$ are examples of f.g.~solvable linear groups.
In \cite{KonigL15} it was shown that for every f.g.~solvable linear group the word problem can be solved in $\TC^0$.
This leads to the question whether for every f.g.~solvable linear group the power word problem belongs to $\TC^0$.
\item The power word problem is a restriction of the compressed word problem, where it is asked whether the word produced by a so-called
straight-line program (a context-free grammar that produces a single word) represents the group identity; see \cite{Loh14}. 
The compressed word problem for $\BS(1,q)$ belongs to $\coRP$ (the complement of randomized polynomial time); this holds in fact
for every f.g.~linear group \cite{Loh14}. No better complexity bound is known for the compressed word problem for $\BS(1,q)$.
\item Is the knapsack problem decidable for every matrix group $T(\alpha)$ with $\alpha \in \mathbb{C}\setminus\{0\}$?
\end{itemize}

\section*{Acknowledgements} Markus Lohrey has been supported by the DFG research project Lo748/12-1.


\begin{thebibliography}{10}

\bibitem{AllenderBD14}
Eric Allender, Nikhil Balaji, and Samir Datta.
\newblock Low-depth uniform threshold circuits and the bit-complexity of
  straight line programs.
\newblock In {\em Proceedings of the 39th International Symposium on
  Mathematical Foundations of Computer Science 2014, {MFCS} 2014}, volume 8635
  of {\em Lecture Notes in Computer Science}, pages 13--24. Springer, 2014.
\newblock \href {https://doi.org/10.1007/978-3-662-44465-8\_2}
  {\path{doi:10.1007/978-3-662-44465-8\_2}}.

\bibitem{AroraBarak09}
Sanjeev Arora and Boaz Barak.
\newblock {\em Computational Complexity - {A} Modern Approach}.
\newblock Cambridge University Press, 2009.
\newblock URL:
  \url{http://www.cambridge.org/catalogue/catalogue.asp?isbn=9780521424264}.

\bibitem{BabaiBCIL96}
L{\'a}zl{\'o} Babai, Robert Beals, Jin yi~Cai, G{\'a}bor Ivanyos, and Eugene
  M.Luks.
\newblock Multiplicative equations over commuting matrices.
\newblock In {\em Proceedings of the Seventh Annual {ACM-SIAM} Symposium on
  Discrete Algorithms, SODA 1996}, pages 498--507. {ACM/SIAM}, 1996.
\newblock URL: \url{http://dl.acm.org/citation.cfm?id=313852.314109}.

\bibitem{BFLW19}
Laurent Bartholdi, Michael Figelius, Markus Lohrey, and Armin Wei{\ss}.
\newblock Groups with {ALOGTIME}-hard word problems and {PSPACE}-complete
  circuit value problems.
\newblock In {\em Proceedings of the 35th Computational Complexity Conference,
  {CCC} 2020}, volume 169 of {\em LIPIcs}, pages 29:1--29:29. Schloss Dagstuhl
  - Leibniz-Zentrum f{\"{u}}r Informatik, 2020.
\newblock \href {https://doi.org/10.4230/LIPIcs.CCC.2020.29}
  {\path{doi:10.4230/LIPIcs.CCC.2020.29}}.

\bibitem{baumslag1961wreath}
Gilbert Baumslag.
\newblock Wreath products and finitely presented groups.
\newblock {\em Mathematische Zeitschrift}, 75(1):22--28, 1961.
\newblock \href {https://doi.org/10.1007/BF01211007}
  {\path{doi:10.1007/BF01211007}}.

\bibitem{BergstrasserGZ21}
Pascal Bergstr{\"{a}}{\ss}er, Moses Ganardi, and Georg Zetzsche.
\newblock A characterization of wreath products where knapsack is decidable.
\newblock In {\em Proceedings of the 38th International Symposium on
  Theoretical Aspects of Computer Science, {STACS} 2021}, volume 187 of {\em
  LIPIcs}, pages 11:1--11:17. Schloss Dagstuhl - Leibniz-Zentrum f{\"{u}}r
  Informatik, 2021.
\newblock \href {https://doi.org/10.4230/LIPIcs.STACS.2021.11}
  {\path{doi:10.4230/LIPIcs.STACS.2021.11}}.

\bibitem{BradfordD88}
Russell~J. Bradford and James~H. Davenport.
\newblock Effective tests for cyclotonic polynomials.
\newblock In {\em Proceedings of the 47th International Symposium on Symbolic
  and Algebraic Computation, ISSAC 1988}, volume 358 of {\em Lecture Notes in
  Computer Science}, pages 244--251. Springer, 1989.
\newblock \href {https://doi.org/10.1007/3-540-51084-2\_22}
  {\path{doi:10.1007/3-540-51084-2\_22}}.

\bibitem{Bruyere1985}
V\'{e}ronique Bruy\`{e}re.
\newblock Entiers et automates finis.
\newblock M\'{e}moire de fin d'\'{e}tudes, Universit\'{e} de Mons, 1985.

\bibitem{bruyere1994logic}
V{\'e}ronique Bruy\`{e}re, Georges Hansel, Christian Michaux, and Roger
  Villemaire.
\newblock Logic and $p$-recognizable sets of integers.
\newblock {\em Bulletin of the Belgian Mathematical Society}, 1:191--238, 1994.
\newblock \href {https://doi.org/10.36045/bbms/1103408547}
  {\path{doi:10.36045/bbms/1103408547}}.

\bibitem{buchi1960weak}
J.~Richard B{\"u}chi.
\newblock Weak second-order arithmetic and finite automata.
\newblock {\em Mathematical Logic Quarterly}, 6(1-6):66--92, 1960.
\newblock \href {https://doi.org/10.1002/malq.19600060105}
  {\path{doi:10.1002/malq.19600060105}}.

\bibitem{BuchiS88}
J.~Richard B{\"{u}}chi and Steven Senger.
\newblock Definability in the existential theory of concatenation and
  undecidable extensions of this theory.
\newblock {\em Mathematical Logic Quarterly}, 34(4):337--342, 1988.
\newblock \href {https://doi.org/10.1002/malq.19880340410}
  {\path{doi:10.1002/malq.19880340410}}.

\bibitem{CaiLZ00}
Jin{-}Yi Cai, Richard~J. Lipton, and Yechezkel Zalcstein.
\newblock The complexity of the {A} {B} {C} problem.
\newblock {\em {SIAM} Journal on Computing}, 29(6):1878--1888, 2000.
\newblock \href {https://doi.org/10.1137/S0097539794276853}
  {\path{doi:10.1137/S0097539794276853}}.

\bibitem{DudkinTreyer2018}
Fedor Dudkin and Alexander Treyer.
\newblock Knapsack problem for {Baumslag--Solitar} groups.
\newblock {\em Siberian Journal of Pure and Applied Mathematics}, 18:43--55,
  2018.
\newblock \href {https://doi.org/10.33048/pam.2018.18.404}
  {\path{doi:10.33048/pam.2018.18.404}}.

\bibitem{Eberly89}
Wayne Eberly.
\newblock Very fast parallel polynomial arithmetic.
\newblock {\em {SIAM} Journal on Computing}, 18(5):955--976, 1989.
\newblock \href {https://doi.org/10.1137/0218066} {\path{doi:10.1137/0218066}}.

\bibitem{FGLZ20}
Michael Figelius, Moses Ganardi, Markus Lohrey, and Georg Zetzsche.
\newblock The complexity of knapsack problems in wreath products.
\newblock In {\em Proceedings of the 47th International Colloquium on Automata,
  Languages, and Programming, {ICALP} 2020}, volume 168 of {\em LIPIcs}, pages
  126:1--126:18. Schloss Dagstuhl - Leibniz-Zentrum f{\"{u}}r Informatik, 2020.
\newblock \href {https://doi.org/10.4230/LIPIcs.ICALP.2020.126}
  {\path{doi:10.4230/LIPIcs.ICALP.2020.126}}.

\bibitem{FiLoZe21}
Michael Figelius, Markus Lohrey, and Georg Zetzsche.
\newblock Closure properties of knapsack semilinear groups.
\newblock {\em Journal of Algebra}, 589(1):437--482, 2022.
\newblock \href {https://doi.org/10.1016/j.jalgebra.2021.08.016}
  {\path{doi:10.1016/j.jalgebra.2021.08.016}}.

\bibitem{FrenkelNU15}
Elizaveta Frenkel, Andrey Nikolaev, and Alexander Ushakov.
\newblock Knapsack problems in products of groups.
\newblock {\em Journal of Symbolic Computation}, 74:96--108, 2016.
\newblock \href {https://doi.org/10.1016/j.jsc.2015.05.006}
  {\path{doi:10.1016/j.jsc.2015.05.006}}.

\bibitem{GalbyOW15}
Esther Galby, Jo{\"{e}}l Ouaknine, and James Worrell.
\newblock On matrix powering in low dimensions.
\newblock In {\em Proceedings of the 32nd International Symposium on
  Theoretical Aspects of Computer Science, {STACS} 2015}, volume~30 of {\em
  LIPIcs}, pages 329--340. Schloss Dagstuhl - Leibniz-Zentrum f{\"{u}}r
  Informatik, 2015.
\newblock \href {https://doi.org/10.4230/LIPIcs.STACS.2015.329}
  {\path{doi:10.4230/LIPIcs.STACS.2015.329}}.

\bibitem{GanardiKLZ18}
Moses Ganardi, Daniel K{\"{o}}nig, Markus Lohrey, and Georg Zetzsche.
\newblock Knapsack problems for wreath products.
\newblock In {\em Proceedings of the 35th Symposium on Theoretical Aspects of
  Computer Science, {STACS} 2018}, volume~96 of {\em LIPIcs}, pages
  32:1--32:13. Schloss Dagstuhl - Leibniz-Zentrum fuer Informatik, 2018.
\newblock \href {https://doi.org/10.4230/LIPIcs.STACS.2018.32}
  {\path{doi:10.4230/LIPIcs.STACS.2018.32}}.

\bibitem{Ge93}
Guoqiang Ge.
\newblock Testing equalities of multiplicative representations in polynomial
  time (extended abstract).
\newblock In {\em Proceedings of the 34th Annual Symposium on Foundations of
  Computer Science, FOCS 1993}, pages 422--426. IEEE Computer Society, 1993.
\newblock \href {https://doi.org/10.1109/SFCS.1993.366845}
  {\path{doi:10.1109/SFCS.1993.366845}}.

\bibitem{GuepinHaaseWorrell2019}
Florent Gu{\'{e}}pin, Christoph Haase, and James Worrell.
\newblock On the existential theories of {B{\"{u}}chi} arithmetic and linear
  $p$-adic fields.
\newblock In {\em Proceedings of the 34th Annual {ACM/IEEE} Symposium on Logic
  in Computer Science, {LICS} 2019}, pages 1--10. IEEE Computer Society, 2019.
\newblock \href {https://doi.org/10.1109/LICS.2019.8785681}
  {\path{doi:10.1109/LICS.2019.8785681}}.

\bibitem{Guyot12}
Luc Guyot.
\newblock Limits of metabelian groups.
\newblock {\em International Journal of Algebra and Computation}, 22(4), 2012.
\newblock \href {https://doi.org/10.1142/S0218196712500312}
  {\path{doi:10.1142/S0218196712500312}}.

\bibitem{Guyot18}
Luc Guyot.
\newblock Generators of split extensions of abelian groups by cyclic groups.
\newblock {\em Groups, Geometry, and Dynamics}, 12(2):765--802, 2018.
\newblock \href {https://doi.org/10.4171/GGD/455} {\path{doi:10.4171/GGD/455}}.

\bibitem{Haa11}
Christoph Haase.
\newblock {\em On the complexity of model checking counter automata}.
\newblock PhD thesis, University of {Oxford}, St Catherine's College, 2011.

\bibitem{HesseAB02}
William Hesse, Eric Allender, and David A.~Mix Barrington.
\newblock Uniform constant-depth threshold circuits for division and iterated
  multiplication.
\newblock {\em Journal of Computer and System Sciences}, 65(4):695--716, 2002.
\newblock \href {https://doi.org/10.1016/S0022-0000(02)00025-9}
  {\path{doi:10.1016/S0022-0000(02)00025-9}}.

\bibitem{Hooley1967}
Christopher Hooley.
\newblock On {Artin}'s conjecture.
\newblock {\em Journal f\"{u}r die reine und angewandte Mathematik},
  1967(225):209--220, 1967.
\newblock \href {https://doi.org/10.1515/crll.1967.225.209}
  {\path{doi:10.1515/crll.1967.225.209}}.

\bibitem{KLM20}
Olga Kharlampovich, Laura L\'{o}pez, and Alexei Myasnikov.
\newblock The {Diophantine} problem in some metabelian groups.
\newblock {\em Mathematics of Computation}, 89:2507--2519, 2020.
\newblock \href {https://doi.org/10.1090/mcom/3533}
  {\path{doi:10.1090/mcom/3533}}.

\bibitem{KonigL15}
Daniel K{\"{o}}nig and Markus Lohrey.
\newblock Evaluation of circuits over nilpotent and polycyclic groups.
\newblock {\em Algorithmica}, 80(5):1459--1492, 2018.
\newblock \href {https://doi.org/10.1007/s00453-017-0343-z}
  {\path{doi:10.1007/s00453-017-0343-z}}.

\bibitem{KoenigLohreyZetzsche2015a}
Daniel K{\"o}nig, Markus Lohrey, and Georg Zetzsche.
\newblock Knapsack and subset sum problems in nilpotent, polycyclic, and
  co-context-free groups.
\newblock In {\em Algebra and Computer Science}, volume 677 of {\em
  Contemporary Mathematics}, pages 138--153. American Mathematical Society,
  2016.
\newblock \href {https://doi.org/10.1090/conm/677/13625}
  {\path{doi:10.1090/conm/677/13625}}.

\bibitem{LaMc98}
Klaus-J{\"o}rn Lange and Pierre McKenzie.
\newblock On the complexity of free monoid morphisms.
\newblock In {\em Proceedings of the 9th International Symposium on Algorithms
  and Computation, {ISAAC} 1998}, number 1533 in Lecture Notes in Computer
  Science, pages 247--256. Springer, 1998.
\newblock \href {https://doi.org/10.1007/3-540-49381-6\_27}
  {\path{doi:10.1007/3-540-49381-6\_27}}.

\bibitem{Lenstra}
H.~W. Lenstra, Jr.
\newblock Finding small degree factors of lacunary polynomials.
\newblock In {\em Number Theory in Progress, vol.~1 Diophantine Problems and
  Polynomials}, pages 267--276. Walter de Gruyter, 1999.

\bibitem{Loh14}
Markus Lohrey.
\newblock {\em The Compressed Word Problem for Groups}.
\newblock SpringerBriefs in Mathematics. Springer, 2014.
\newblock \href {https://doi.org/10.1007/978-1-4939-0748-9}
  {\path{doi:10.1007/978-1-4939-0748-9}}.

\bibitem{Loh19hyp}
Markus Lohrey.
\newblock Knapsack in hyperbolic groups.
\newblock {\em Journal of Algebra}, 545:390--415, 2020.
\newblock \href {https://doi.org/10.1016/j.jalgebra.2019.04.008}
  {\path{doi:10.1016/j.jalgebra.2019.04.008}}.

\bibitem{LoWe19}
Markus Lohrey and Armin Wei{\ss}.
\newblock The power word problem.
\newblock In {\em Proceedings of the 44th International Symposium on
  Mathematical Foundations of Computer Science, MFCS 2019}, volume 138 of {\em
  LIPIcs}, pages 43:1--43:15. Schloss Dagstuhl - Leibniz-Zentrum f{\"u}r
  Informatik, 2019.
\newblock \href {https://doi.org/10.4230/LIPIcs.MFCS.2019.43}
  {\path{doi:10.4230/LIPIcs.MFCS.2019.43}}.

\bibitem{LohreyZ18}
Markus Lohrey and Georg Zetzsche.
\newblock Knapsack in graph groups.
\newblock {\em Theory of Computing Systems}, 62(1):192--246, 2018.
\newblock \href {https://doi.org/10.1007/s00224-017-9808-3}
  {\path{doi:10.1007/s00224-017-9808-3}}.

\bibitem{LohreyZ20}
Markus Lohrey and Georg Zetzsche.
\newblock Knapsack and the power word problem in solvable {Baumslag-Solitar}
  groups.
\newblock In {\em Proceedings of the 45th International Symposium on
  Mathematical Foundations of Computer Science, {MFCS} 2020}, volume 170 of
  {\em LIPIcs}, pages 67:1--67:15. Schloss Dagstuhl - Leibniz-Zentrum f{\"{u}}r
  Informatik, 2020.
\newblock \href {https://doi.org/10.4230/LIPIcs.MFCS.2020.67}
  {\path{doi:10.4230/LIPIcs.MFCS.2020.67}}.

\bibitem{MRUV10}
Alexei Miasnikov, Vitaly Roman'kov, Alexander Ushakov, and Anatoly Vershik.
\newblock The word and geodesic problems in free solvable groups.
\newblock {\em Transactions of the American Mathematical Society},
  362(9):4655--4682, 2010.
\newblock \href {https://doi.org/10.1090/S0002-9947-10-04959-7}
  {\path{doi:10.1090/S0002-9947-10-04959-7}}.

\bibitem{MishchenkoT17}
Alexei Mishchenko and Alexander Treier.
\newblock Knapsack problem for nilpotent groups.
\newblock {\em Groups Complexity Cryptology}, 9(1):87, 2017.
\newblock \href {https://doi.org/10.1515/gcc-2017-0006}
  {\path{doi:10.1515/gcc-2017-0006}}.

\bibitem{MyNiUs14}
Alexei Myasnikov, Andrey Nikolaev, and Alexander Ushakov.
\newblock Knapsack problems in groups.
\newblock {\em Mathematics of Computation}, 84:987--1016, 2015.
\newblock \href {https://doi.org/10.1090/S0025-5718-2014-02880-9}
  {\path{doi:10.1090/S0025-5718-2014-02880-9}}.

\bibitem{Nathanson2000}
Melvyn~B. Nathanson.
\newblock {\em Elementary Methods in Number Theory}.
\newblock Springer, 2000.
\newblock \href {https://doi.org/10.1007/b98870} {\path{doi:10.1007/b98870}}.

\bibitem{StoberW22}
Florian Stober and Armin Wei{\ss}.
\newblock The power word problem in graph products.
\newblock In {\em Proceedings to the 26th International Conference on
  Developments in Language Theory, {DLT} 2022}, volume 13257 of {\em Lecture
  Notes in Computer Science}, pages 286--298. Springer, 2022.
\newblock \href {https://doi.org/10.1007/978-3-031-05578-2\_23}
  {\path{doi:10.1007/978-3-031-05578-2\_23}}.

\bibitem{Vol99}
Heribert Vollmer.
\newblock {\em Introduction to Circuit Complexity}.
\newblock Springer, 1999.
\newblock \href {https://doi.org/10.1007/978-3-662-03927-4}
  {\path{doi:10.1007/978-3-662-03927-4}}.

\bibitem{Weiss15b}
Armin Wei{\ss}.
\newblock {\em On the complexity of conjugacy in amalgamated products and {HNN}
  extensions}.
\newblock PhD thesis, University of Stuttgart, 2015.
\newblock URL: \url{http://elib.uni-stuttgart.de/opus/volltexte/2015/10018/}.

\bibitem{Weiss16}
Armin Wei{\ss}.
\newblock A logspace solution to the word and conjugacy problem of generalized
  {Baumslag-Solitar} groups.
\newblock In {\em Algebra and Computer Science}, volume 677 of {\em
  Contemporary Mathematics}. American Mathematical Society, 2016.
\newblock \href {https://doi.org/https://doi.org/10.1090/conm/677/13628}
  {\path{doi:https://doi.org/10.1090/conm/677/13628}}.

\bibitem{Woe00}
Wolfgang Woess.
\newblock {\em Random Walks on Infinite Graphs and Groups}.
\newblock Cambridge University Press, 2000.
\newblock \href {https://doi.org/10.1017/CBO9780511470967}
  {\path{doi:10.1017/CBO9780511470967}}.

\end{thebibliography}

 \end{document}